\documentclass[a4paper,10pt]{amsart}

\usepackage[all]{xy}
\usepackage[T1]{fontenc}
\usepackage[leqno]{amsmath}
\usepackage{amssymb}
\usepackage{amsthm}
\usepackage[english]{babel}

\newcounter{eq}[subsection]
\newcounter{seq}[subsubsection]

\newcommand\seq{\stepcounter{seq}}
\newcommand\snum{\seq\tag{\thesubsubsection.\theseq}}

\theoremstyle{plain}
\newtheorem*{theo}{Theorem}

\newtheorem{lemma}[subsection]{Lemma}
\newtheorem{prop}[subsection]{Proposition}

\theoremstyle{definition}
\newtheorem{rem}[subsection]{Remark}

\theoremstyle{plain}
\newtheorem{stheoreme}[subsubsection]{Theorem}
\newtheorem{slemma}[subsubsection]{Lemma}
\newtheorem{sprop}[subsubsection]{Proposition}
\newtheorem{scor}[subsubsection]{Corollary}
\theoremstyle{definition}
\newtheorem{srem}[subsubsection]{Remark}

\newcommand{\transp}[1]{{}^{\mathit t}\!{#1}}

\newcommand\lra{\longrightarrow}
\newcommand\mso{\mathcal{M}_{\mathbf{SO}_r}}
\newcommand\mo{\mathcal{M}_{\mathbf{O}_r}}
\newcommand\msl{\mathcal{M}_{\mathbf{SL}_r}}
\newcommand\sor{{\mathbf{SO}_r}}
\newcommand\slr{{\mathbf{SL}_r}}
\newcommand\mgl{\mathcal{M}_{\mathbf{GL}_r}}
\newcommand\glr{{\mathbf{GL}_r}}
\newcommand\gl{{\mathbf{GL}}}
\newcommand\ssl{{\mathbf{SL}}}
\renewcommand\o{{\mathbf{O}}}
\newcommand\so{{\mathbf{SO}}}
\newcommand\Sp{{\mathbf{Sp}}}
\newcommand{\sok}{{\boldsymbol{\mathfrak{so}}}}
\newcommand{\glk}{{\boldsymbol{\mathfrak{gl}}}}
\newcommand{\spk}{{\boldsymbol{\mathfrak{sp}}}}
\renewcommand\O{\mathcal{O}}

\newcommand\git{/\!\!/}

\newcommand\Q{\widetilde Q}
\newcommand\A{\tilde A}
\newcommand\X{\widetilde X}
\newcommand\St{\tilde S}
\newcommand\R{\widetilde R}

\newcommand\G{{\o_N \times \Sp_{N'}}}
\newcommand\M{\mathrm M}
\newcommand\tr{\mathrm{tr}}
\newcommand\n{{N+N'}}

\author{Olivier Serman}
\address{Laboratoire J.-A. Dieudonn\'e \\ UMR 6621 du CNRS \\ Universit\'e de Nice \\ Parc Valrose \\ F-06108 Nice Cedex 02}
\email{serman@math.unice.fr}
\title{Moduli spaces of orthogonal bundles over an algebraic curve}

\begin{document}

\begin{abstract}
We prove that, given $C$ a smooth projective curve of genus $g \geqslant 2$, the forgetful map $\mo \to \mgl$ from the moduli space of orthogonal bundles to the moduli space of all vector bundles over $C$ is an embedding. Our proof relies on an explicit description of a set of generators for the polynomial invariants on the representation space of a quiver under the action of a product of classical groups. 
\end{abstract}

\maketitle

\section*{Introduction}

Let $C$ be a smooth, irreducible, projective algebraic curve of genus $g \geqslant 2$ over an algebraically closed field $k$ of characteristic zero. If $G$ is a reductive group we denote by $\mathcal{M}_G$ the moduli space of semi-stable principal $G$-bundles on $C$. 

We focus here on the case $G=\sor$, which amounts to considering the moduli space of semi-stable orthogonal bundles of rank $r$ with an orientation. It is a normal projective variety, composed of two connected components distinguished by the second Stiefel-Whitney class. This space is related to the moduli space $\msl$ of vector bundles of rank $r$ and trivial determinant on $C$ through the forgetful map $\mso \to \msl$ which sends any $\sor$-bundle to its underlying vector bundle. It is natural to ask whether this map is a closed embedding. In fact, when $r$ is even, it even fails to be injective, and it is therefore more convenient to ask the same question about $\mo \to \mgl$. 

In the same way we consider the forgetful morphism $\mathcal M_{{\Sp}_{2r}} \to \mathcal M_{\ssl_{2r}}$ defined on the moduli space of symplectic bundles of rank $2r$ on $C$.

Our main result may be stated as follows:

\begin{theo}
{\rm{(i)}} The forgetful map $\mo \rightarrow \mgl$ is an embedding. 

\indent {\rm {(ii)}} When $r$ is odd, $\mso \to \msl$ is again an embedding, while, when $r$ is even, it is a $2$-sheeted cover onto its image.

\indent {\rm{(iii)}} The forgetful map $\mathcal M_{{\Sp}_{2r}} \to \mathcal M_{\ssl_{2r}}$ is also an embedding.
\end{theo}

We give the full proof for the orthogonal case, and sketch the obvious modifications required by the symplectic one.

We consider in the first section the injectivity of $\mso \to \msl$: this comes down to a close comparison of the equivalence relations between $\sor$-bundles and vector bundles which define the closed points of the corresponding moduli spaces. We then check that the tangent maps of $\mo \to \mgl$ are injective. This differential point of view is much more involved: it relies on Luna's slice \'etale theorem, which naturally leads to the consideration of representations of quivers. To carry our discussion to its end we need an auxiliary result relative to the invariant theory of these representations for the action of a product of classical groups: this is the aim of the second section (note that a characteristic-free proof of this result can be found in \cite{lopatin}). In the third one we show first how this computation results in our main theorem. We then give a few complements about the local structure of the moduli space $\mo$, in the same way as Laszlo did with $\mgl$ (see \cite{localstr}).

\section{Injectivity of $\mathcal M _{\sor} \rightarrow \mathcal M _{\slr}$}

In this section we study the injectivity of the forgetful map $\mathcal M _{\sor} \lra \mathcal M _{\slr}$, which is already known to be finite (e.g. by \cite[8.5]{balsesh}). The closed points of $\mathcal M _G$ are in a one-to-one correspondence with the set of equivalence classes of semi-stable $G$-bundles (cf. \cite{ramanathanthesis}). When $G=\slr$ one easily recovers from this notion Seshadri's definition of $S$-equivalence for vector bundles. We proceed here in the same way to link together equivalence between $\sor$-bundles and $S$-equivalence between their underlying vector bundles.   

\subsection{} Let $Q$ be a quadratic form on a $r$-dimensional vector space $V$. Any parabolic subgroup of $\so(Q)$ is the stabilizer of an isotropic flag of $V$. If $0=N_0 \subset N_1 \subset \cdots \subset N_l$ is such a flag then its stabilizer $\Gamma$ stabilizes $0=N_0 \subset N_1 \subset \cdots \subset N_l \subset N_l^\bot \subset \cdots \subset N_1^\bot \subset V$ too. In a basis adapted to the first flag, in which $Q$ is represented by the matrix $\begin{pmatrix} 0 & 0 & \mathrm I_{\dim N_l} \\ 0 & Q' & 0 \\ \mathrm I_{\dim N_l} & 0 & 0 \end{pmatrix}$,  $\Gamma$ must be a subgroup of the group of matrices of the form
$${M=\left( \begin{array}{ccc}
\begin{array}{cccc}{A_1} & \ast  &\cdots  &\ast   \\ 0     & { A_2}&\ddots  &\vdots  \\ \vdots&\ddots & \ddots & \ast  \\ 0  & \cdots&  0 &  A_l \end{array} & \begin{array}{c} \ast \\ \vdots \\ \vdots\\ \ast \end{array} &\textrm{\huge $\ast$} \\
\begin{array}{cccc} 0\ &\cdots \ &\cdots \!&\ 0 \end{array}  & B & \begin{array}{cccc} \ast \quad  & \ \cdots\!  &\ \quad \cdots &\quad \ast  \end{array} \\
\textrm{\LARGE  $0$}  & \begin{array}{c} 0\\ \vdots\\ \vdots\\0\end{array} & \begin{array}{cccc} \ \transp{A_1}\! ^{-1} & 0  &\cdots& 0 \\  \ast &\ \transp{A_2}\! ^{-1}&\ddots& \vdots    \\ \vdots & \ddots & \ddots& 0  \\ \ast & \cdots & \ast &\ \transp{A_l}\! ^{-1} \end{array}
\end{array}\right)},$$
\noindent where $A_i \in \mathbf{GL}_{r_i - r_{i-1}}$, $B \in \mathbf{SO}(Q')$ (and $r_i=\dim N_i$).

Hence, if $P$ is a semi-stable $\sor$-bundle and $E=P(\slr)$ its associated vector bundle, a reduction of structure group $\sigma$ of $P$ to a parabolic subgroup $\Gamma$ defined by an isotropic flag  $0=N_0 \subset N_1 \subset \cdots \subset N_l$ induces isotropic subbundles $E_i = \sigma^\ast E (N_i) \subset E$ of rank $r_i$ giving a filtration  $0=E_0 \subset E_1 \subset \cdots \subset E_l \subset E_l^\bot \subset \cdots \subset E_1^\bot \subset E$. But one knows how to construct such a filtration, with the extra conditions that $E_i/E_{i-1}$ is a stable bundle of degree $0$ and  $E_l ^\bot / E_l$ a stable orthogonal bundle. By \cite[4.5]{ramanan} the latter splits as a direct orthogonal sum of mutually non-isomorphic stable bundles. The graded object $\mathrm{gr} E_{\bullet}$ is then precisely the Jordan-H\"older one, which is known to characterize the point of $\msl$ corresponding to $E$. \label{jordan-holder} 

\subsection{} A representative of the equivalence class of $P$ as an orthogonal bundle is given by the $\sor$-bundle $\mathrm{gr} P$ obtained from a suitable reduction of structure group of $P$ to a parabolic subgroup of $\sor$ (cf. \cite[3.12]{ramanathanthesis}). Let us check that the reduction $\sigma$ attached to the above filtration satisfies the conditions of (\textit{loc. cit.}). As any character $\chi$ on $\Gamma$ which is trivial on the neutral component $Z^0$ of its center is of the form $M \mapsto \prod \det(A_i)^{\alpha_i}$ with $\sum (\dim N_i / N_{i-1}) \alpha_i=0$, the assertion $\deg(\chi_\ast \sigma^\ast E (k))=0$ is equivalent to $\mu(E_i / E_{i-1})=0$, $i=1,\ldots,l$, and $\sigma$ definitely is an admissible reduction. 
 
On the other hand the unipotent radical  $\mathrm{R}_u(\Gamma)$ of $\Gamma$ is the unipotent part of the (neutral component of the) intersection of all of its Borel subgroups. Since these are the stabilizers of the flags adapted to the one giving $\Gamma$ whose length is maximal, $\mathrm{R}_u(\Gamma)$ is a subgroup consisting of matrices $M$ with $A_i=\mathrm{id}$ for all $i$ and $B=\mathrm{id}$. We deduce from the preceding a Levi decomposition of $\Gamma$, the Levi component $L$ being isomorphic to  the product  $\prod \mathbf{GL}_{r_i-r_{i-1}} \times \mathbf{SO}_{r-2r_l}$. Let $p \colon \Gamma \rightarrow L$ be the projection on this Levi component. The required stability of the $L$-bundle  $p_\ast \sigma^\ast E$ is then a consequence of the stability of the successive quotients $E_i / E_{i-1}$ together with the following easy fact:

\begin{lemma}
Let $(E_i)_{i=1,2}$ be two stable $G_i$-bundles. Then the  $G_1\times G_2$-bundle $E_1 \times_C E_2$ is also stable.
\end{lemma}

According to \cite[11.14\! (1)]{borel} any parabolic subgroup of $G_1 \times G_2$ is a product  $\Gamma=\Gamma_1 \times \Gamma_2$ where $\Gamma_i$ is a parabolic subgroup of $G_i$. It follows that maximal proper parabolic subgroups may be written $\Gamma_1 \times G_2$ or $G_1 \times \Gamma_2$, $\Gamma_i$ being any maximal proper parabolic subgroup of $G_i$. Therefore the associated bundle $E/\Gamma$ is isomorphic to $E_1/\Gamma_1$ or $E_2/\Gamma_2$, and the lemma is a consequence of the very definition of stability.

\subsection{} The class of $P$ is therefore defined by the $\sor$-bundle $\mathrm{gr} P = p_\ast \sigma^\ast P (\sor)$, and $(\mathrm{gr} P)(\slr)$ is again the Jordan-H\"older graded object. Two semi-stable $\sor$-bundles $P$ and $P'$ are sent to the same point of $\msl$ if and only if $\mathrm{gr} P$ and $\mathrm{gr} P'$ are both obtained from reduction of structure group to $\sor$ of the same polystable vector bundle $E$. Such a reduction amounts to a section of $E/\sor \to C$, and two of them give the same bundles if and only if they are conjugated by the action of  $\mathrm{Aut}_{\slr}(E)$ on $\Gamma(C , E / \sor)$. Elements of $\Gamma(C , E / \sor)$ correspond to isomorphisms $\iota \colon E \buildrel\sim\over\longrightarrow E^\ast$ such that $\iota^\ast=\iota$ and $\det \iota$ is the square of the trivialisation of $\det E$ inherited from the $\slr$-torsor structure. The action of $\mathrm{Aut}_{\slr}(E)$ simply is $$(f , \iota) \in \mathrm{Aut}_{\slr}(E)\times\Gamma(C , E / \sor) \mapsto   f^\ast \iota f.$$ 

\label{pointsfermes} Since $E$ is polystable $\mathrm{Aut}_{\glr}(E)$ acts transitively on the set $\Gamma(C,E/\o_r)$ of all symmetric isomorphisms from $E$ onto $E^\ast$ : indeed the Jordan-H\"older fitration allows us to split $E$ as     
$$E=\bigoplus\limits_i \Bigl(F^{(1)}_{i}\otimes V^{(1)}_{i}\Bigr) \oplus \bigoplus\limits_j \Bigl(F^{(2)}_{j} \otimes V^{(2)}_{j}\Bigr) \oplus \bigoplus\limits_k \Bigl((F^{(3)}_{k} \oplus {F^{(3)}_{k}}{}^\ast)\otimes V^{(3)}_{k}\Bigr),$$
\noindent where $V^{(l)}_{i}$ are finite-dimensional vector spaces and the $F^{(1)}_{i}$ (resp. $F^{(2)}_{j}$, resp. $F^{(3)}_{k}$) are orthogonal (resp. symplectic, resp. non isomorphic to their dual) mutually non isomorphic stable vector bundles, in such a way that any symmetric isomorphism $E \to E^\ast$ is equivalent to the data of orthogonal (resp. symplectic, resp. non degenerate) forms on each one of the $V^{(1)}_i$ (resp. $V^{(2)}_j$, resp. $V^{(3)}_k$). The action of 
$$\mathrm{Aut}_\glr(E)= \prod \gl(V^{(1)}_{i}) \times \prod \gl(V^{(2)}_{j}) \times \bigl(\prod \gl(V^{(3)}_{k}) \times \gl(V^{(3)}_{k}) \Bigr)$$ 
\noindent on the set of these collections is obviously transitive.

Any two elements $\sigma$ and $\sigma'$ of $\Gamma(C , E / \sor)$ are then conjugate under the action of  $\mathrm{Aut}_\glr(E)$, by an automorphism whose determinant equals to $\pm 1$. When $r$ is odd $-\mathrm{id}_E$ is an $\o_r$-isomorphism which exchanges orientation, and the action of $\mathrm{Aut}_{\slr}(E)$ on $\Gamma(C , E / \sor)$ is still transitive. On the contrary when $r$ is even this action fails to remain transitive. For example let $F$ be a vector bundle of rank $r/2$, non isomorphic to its dual, and consider the two orthogonal bundles $F \oplus F ^\ast$ and $F ^\ast \oplus F$, equipped with the standard hyperbolic pairing: these bundles cannot be $\sor$-isomorphic (in fact any orthogonal isomorphism must exchange the orientation). We have proven so far: 

\begin{prop}
When $r$ is odd the map $\mathcal M _{\sor}(k) \rightarrow \mathcal M _{\slr}(k)$ is injective; when $r$ is even this is a finite map of degree $2$.
\end{prop}

\begin{rem} \label{giraud}
The distinction between the odd or even case relies on the fact that the semi-direct product $1 \to \sor \to \mathbf{O}_r \to \mathbb Z / 2 \mathbb Z \to 0$ may be direct or not. When it is direct the map  $H^1(C,\sor) \to H^1(C,\mathbf{O}_r)$ is injective. But as soon as $r$ is even the section is no longer compatible with the action of $\sor$ on the previous exact sequence by inner automorphisms; we know then how to compute the fibers of  $H^1(C,\sor) \to H^1(C,\mathbf{O}_r)$ by twisting this sequence (cf. \cite[III 3.3.4]{giraud}). If $r$ is even we have just chosen a bundle $E$ which induces a non trivial connecting homomorphism (since $\mathrm{Aut}_{\sor}(E)\to \mathrm{Aut}_{\mathbf O _r}(E)$ is then an isomorphism), whence the lack of injectivity.
\end{rem}

\section{Invariant theory of representations of quivers} 

We refer the reader to \cite{LBP} about the notion of representations of a quiver $Q$ of given dimension $\alpha \in \mathbb N^n$ and their isomorphism classes. Let $Q$ stand for a quiver consisting of $n=n_1+n_2+2n_3$ vertices 
$$s_1, \ldots, s_{n_1}, t_1,\ldots,t_{n_2}, u_1,u_1^\ast,\ldots,u_{n_3},u_{n_3}^\ast,$$ 
\noindent and $\alpha \in \mathbb N^n$ be an admissible dimension vector (that is a vector such that  $\alpha_{t_j}$ is even and $\alpha_{u_k}=\alpha_{u_k^\ast}$). We define $\Gamma_\alpha$ to be the group 
$$\Gamma_\alpha=\prod \o_{\alpha_{s_i}} \times \prod \Sp_{\alpha_{t_j}} \times \prod \gl_{\alpha_{u_k}},$$ 
\noindent which is actually thought here as a subgroup of $\gl(\alpha)=\prod_{i=1}^n \gl_{\alpha_i}$ via the inclusions $P \in \gl_{\alpha_{u_k}} \mapsto (P,\transp P^{-1}) \in \gl_{\alpha_{u_k}} \times \gl_{\alpha_{u_k^\ast}}$. The natural action of $\gl(\alpha)$ on the space $R(Q,\alpha)$ of all representations of $Q$ of dimension $\alpha$ restricts to an action of $\Gamma_\alpha$ on $R(Q,\alpha)$. 

Le Bruyn and Procesi have shown in (\textit{loc. cit.}) that the algebra of polynomial invariants  $k[R(Q,\alpha)]^{\gl(\alpha)}$ is generated by traces along oriented cycles in the quiver $Q$. Following their proof, we produce here a set of generators for the algebra of invariants under the action of $\Gamma_\alpha$. The local study of the map $\mo \to \mgl$ made later heavily rests on this description.

\subsection{First fundamental theorem for $\G$}
In this section we adapt the argument of \cite[Appendix 1]{abp}. We will denote by $\M_n$ the space of $n \times n$ matrices. Let $M$ be the matrix  $\left({I_N \atop 0}{0 \atop J_{N'}}\right)$, with $J = \left({0 \atop -I}{I \atop 0}\right)$; it represents a bilinear pairing, given as the standard orthogonal sum of a quadratic form of rank $N$ and a symplectic form of rank $N'$. The key lemma becomes (note that $\M_N \times \M_{N'}$ is identified with its image in $\M_\n$): 

\begin{slemma}
Any polynomial function $f \colon (\M_N \times \M_{N'})(k) \to k$ such that $f(BA)=f(A)$ for all $B \in \G$ may be written $f \colon A \mapsto F(\transp AMA)$ with $F$ a polynomial map on $(\M_N \times \M_{N'})(k) \subset \M_\n(k)$.
\end{slemma}

In other words $f$ factors through the application
$$\pi \colon A \in (\M_N \times \M_{N'})(k) \mapsto \transp AMA \in (\M_N \times \M_{N'})(k).$$ 
\noindent Let $\Psi_{N,N'}$ be its image, which is nothing else than the product of the space of symmetric $N \times N$ matrices and the space of antisymmetric $N' \times N'$ matrices. The restriction of $\pi$ to   $\gl_N \times \gl_{N'}$ identifies the open subset $\Psi_{N,N'}^{\circ}$ consisting of non-degenerate forms with the geometric quotient $(\gl_N \times \gl_{N'}) \git (\G)$ (say by \cite[Proposition 0.2]{GIT}). The lemma then follows from the commutative diagram 
$$\xymatrix@C=50pt{
(\M_N \times \M_{N'}) \git (\G) \ar[r]_-{\pi} & \Psi_{N,N'} \ {\vphantom{\ ;}}\\
(\gl_N \times \gl_{N'}) \git (\G) \ar@{}[u]|{\bigcup} \ar[r]^-\sim_-\pi & {\vphantom{}}\  \Psi_{N,N'}^{\circ} \,;  \ar@{}[u]|{\bigcup}
}$$
\noindent the restriction to $(\gl_N \times \gl_{N'}) \git (\G)$ of a map $f$ defined on the good quotient $(\M_N \times \M_{N'}) \git (\G)$ must indeed be induced by a function of the form $A \in \gl_N \times \gl_{N'} \mapsto F(\transp AMA)/H(\transp AM A)$ with $F$ and $H$ two coprime polynomials (defined on $(\M_N \times \M_{N'})(k)$). The equality $F(\transp AMA)=f(A)H(\transp AMA)$ finally ensures that $H$ is invertible.

\subsubsection{} \label{firstmaintheorem} We are now in a position to establish the first main theorem of invariant theory for $\G$. Let us start with the case of multilinear invariants, again after \cite{abp}. Let $V$ be a vector space of dimension $\n$ endowed with the non-degenerate bilinear form  $\langle\cdot\,,\cdot\rangle$ given by the matrix $M$. So $V = V_1 \buildrel\bot\over\oplus V_2$, $V_1$ being a quadratic space of dimension $N$ and $V_2$ a symplectic space of dimension $N'$.

\begin{stheoreme} \label{firsttheorem}
Any linear $\G$-invariant morphism $V^{\otimes 2i} \to k$ is a linear combination of functions of the form
$$\varphi_\sigma \colon v_1 \otimes \cdots \otimes v_{2i} \longmapsto \langle v_{\sigma(1)},v_{\sigma(2)} \rangle \cdots \langle v_{\sigma(2i-1)},v_{\sigma(2i)} \rangle ,$$
with $\sigma \in \mathfrak S _{2i}$.
\end{stheoreme}

Let $\varphi \colon V^{\otimes 2i} \to k$ be any linear $\G$-invariant map, and consider the following polynomial function: 
$$f \colon (A,\omega) \in (\mathrm{End}\, V_1 \oplus \mathrm{End}\, V_2) \times V^{\otimes 2i} \mapsto \varphi(A\omega) \in k.$$ 
By the previous lemma there exists a polynomial $F$ on $(\mathbf S^2 V_1^\ast \oplus \boldsymbol \Lambda^2 V_2^\ast) \times V^{\otimes 2i}$, linear in the second variable, such that $f(A,\omega)= F(\transp A M A,\omega)$. This polynomial certainly is invariant for the natural action of $\gl (V_1) \times \gl (V_2)$ on $(\mathbf S^2 V_1^\ast \oplus \boldsymbol \Lambda^2 V_2^\ast) \times V^{\otimes 2i}$: for any $\Gamma \in \gl(V_1) \times \gl(V_2)$, we have $F({}^{\mathit{t}}\Gamma^{-1} \transp A M A \Gamma^{-1}, \Gamma \omega)=F(\transp A M A,\omega)$. 

The assertion results, by polarization, from the description of linear forms on ${V_1^\ast}^{\otimes a_1} \otimes V_1^{\otimes b_1} \otimes {V_2^\ast}^{\otimes a_2} \otimes V_2^{\otimes b_2}$ which are invariant for the action of  ${\gl(V_1)\times\gl(V_2)}$: $F$ is an homogeneous function of degree $i$ in its first variable, which therefore arises from linear forms on $(\mathbf S^2V_1^\ast)^{\otimes k} \otimes V_1^{\otimes 2k} \otimes (\boldsymbol\Lambda^2V_2^\ast)^{\otimes i-k}\otimes V_2^{\otimes 2i-2k}$ (via the projections $(\mathbf S^2 V_1^\ast \oplus \boldsymbol \Lambda^2 V_2^\ast)^{\otimes i} \times V^{\otimes 2i} \to (\mathbf S^2V_1^\ast)^{\otimes k} \otimes V_1^{\otimes 2k} \otimes (\boldsymbol\Lambda^2V_2^\ast)^{\otimes i-k}\otimes V_2^{\otimes 2i-2k}$). Since $\varphi(\omega)=F(M,\omega)$, one then just has to evaluate $F$ on $M$. 

\subsubsection{} \label{procesitrick} One easily deduces from the foregoing a family of generators for the algebra of polynomial invariants under the diagonal action (by conjugation) of $\G$ on $\M_{\n}(k)^i$: according to \cite[\S7]{PrAdv} it is enough to work out the behaviour of the composition, the trace and the adjunction (denoted by  $A \mapsto A^\ast=M^{-1}\transp A M$) via the identification $\mathrm{End}\,V \simeq V \otimes V$ induced by the bilinear pairing. If $v=v_1+v_2 \in V=V_1 \oplus V_2$ (cf. \ref{firstmaintheorem}), we have the following identities:

\begin{itemize}
\item $(v \otimes w) \circ (u \otimes t)=\langle v,t\rangle u \otimes w$,
\item $\tr(v \otimes w)=\langle v,w \rangle$,
\item $((v_1+v_2) \otimes w)^\ast = w \otimes (v_1-v_2)$.
\end{itemize}

This relations allow us to translate the functions $\varphi_\sigma$ occuring in theorem \ref{firsttheorem} in a way leading to the following statement:

\begin{stheoreme} \label{maintheorem}
Any $\G$-invariant function defined on $\M_{\n}(k)^i$ is a polynomial in the 
$$(A_1,\ldots,A_i) \mapsto \tr(U_{j_1} U_{j_2} \cdots U_{j_k}),$$
with $U_j \in \{A_j,{A_j}^\ast\}$.

The ring of $\G$-equivariant morphisms from $\M_{\n}(k)^i$ to $\M_{\n}(k)$ is generated as an algebra over $k[\M_{\n}(k)^i]^{\G}$ by the elements $(A_1,\ldots,A_i) \mapsto A_j$ or $(A_1,\ldots,A_i) \mapsto {A_j}^\ast$.
\end{stheoreme}

(The second assertion is implied by the first exactly as in \cite{PrAdv}.)

\subsection{Generalization of the main result of \cite{PrJofalg}} 
The main result of \cite{PrJofalg} asserts that, if $R$ is a $k$-algebra with trace satisfying the $n$-th Cayley-Hamilton identities, then there exists a universal map $R \to \M_{n}(A)$ which induces an isomorphism $R \to \M_n(A)^{\gl_n}$. In this section we state a similar result dealing with algebras with a trace and an antimorphism of order dividing $4$.  

\subsubsection{} \label{foncteur} Let $R$ be a $k$-algebra endowed with an antimorphism $\tau$ of order $4$ (from now on we will write ``of order $4$'' instead of ``of order dividing $4$''). When $R$ is the algebra $\M_\n(B)$ of all matrices with coefficients in a commutative ring $B$ we choose $\tau$ to be the adjunction map (for the considered bilinear form)  $\iota \colon A \in \M_\n(k) \mapsto M^{-1}\transp A M$. As soon as $N$ or $N'$ is zero $\iota$ is in fact of order $2$, and we could restrict ourselves to algebras with anti-involution.

Let $j \colon R \to \M_\n(A)$ be the universal map corresponding to the functor 
$$X_{R,\n} \colon B \mapsto \mathrm{Hom}_k(R,\M_\n(B))$$ 
\noindent (cf. \cite[\S3]{LBP}). The existence of this morphism easily leads to the representability of the functor $\X_{R,\n}$ which associates to any commutative algebra $B$ the set of all morphisms of $k$-algebras from $R$ to $\M_\n(B)$ preserving the antimorphisms. This functor is actually represented by a closed subscheme of $X_{R,\n}$, still called $\X_{R,\n}$: the map $r \in R \mapsto \iota{j\tau^3(r)} \in \M_\n(A)$ comes from a morphism $t \colon A \to A$ of order $4$, and the induced map $\tilde j \colon R \to \M_\n(\A)$ (where $\A$ is the quotient of $A$ by the action of $t$) is universal among the morphisms $R \to \M_\n(B)$ commuting with $\tau$ and $\iota$. 

The group $\G$ acts by conjugation on $\M_\n(B)$, inducing a right action on $\A$, hence an action of $\G$ on $\M_\n(\A)$. The universal map $\tilde j$ maps $R$ to the algebra $\M_\n(\A)^{\G}$ of $\G$-equivariant morphisms from $\X_{R,\n}$ to $\M_\n(k)$ (cf. \cite[1.2]{PrJofalg}).

\subsubsection{} A $k$-algebra with trace and antimorphism of order $4$ is an algebra with trace (cf. \cite{PrJofalg}) equipped with an antimorphism $\tau \colon R \to R$ of order $4$ commuting with $\tr$. The algebra $\M_\n(B)$ carries its usual trace and the antimorphism $\iota$ described earlier. 

Our purpose is to generalize the main theorem of (\textit{loc. cit.}) to any algebra which is a quotient of the algebra $T_{N,N'}$ of $\G$-equivariant morphisms $\M_{\n}(k)^i \to \M_{\n}(k)$. Note that a set of generators for $T_{N,N'}$ has been given in \ref{maintheorem}. 

\begin{sprop}
Let $R$ be a $k$-algebra with trace and antimorphism of order $4$. If $R$ is a quotient of $T_{N,N'}$ then the universal map $\tilde j \colon R \to \M_\n(\A)^{\G}$ is an isomorphism. 
\end{sprop}

It follows from an immediate adaptation of Procesi's proof.  

\subsection{Generators for $k[R(Q,\alpha)]^{\Gamma_\alpha}$}

\subsubsection{} \label{associatequiver} Let us go back to the quiver $Q$, and the action of $\Gamma_\alpha$ on its representation space. Consider the quiver $\Q$ obtained from $Q$ by adding one new arrow $a^\ast \colon \sigma(v') \to \sigma(v)$ for any arrow  $a \colon v \to v'$, where we called $\sigma$ the involution of the set of vertices fixing the $s_i$ and $t_j$, and permuting $u_k$ and $u_k^\ast$. Let $R$ (resp. $\R$) be the path algebra of the opposite quiver $Q^{\mathrm{op}}$ (resp. $\Q^{\mathrm{op}}$). There is quite a natural way to endow $\R$ with an antimorphism $\tau$ of order $4$ whose action on the idempotents (who are associated to the constant paths) comes from the one of $\sigma$: $\tau$ fixes the $e_{s_i}$ and $e_{t_j}$, permutes $e_{u_k}$ and $e_{u_k^\ast}$, while it sends an arrow $a$ to $\varepsilon a^\ast$, where $\varepsilon$ equals $-1$ if $a$ starts from $s_i$, $u_k$ or $u_k^\ast$ and ends at $t_j$, and $1$ otherwise.   

We need here to adapt the map $\iota$ previously defined on $\M_\n(k)$: $\iota$ still associates to a map its adjoint, but the bilinear pairing has to be replaced by the one represented by a matrix of the form 
$$\Phi=\left(\begin{array}{ccc} I_{N_1} &0 &0\\0& J_{N_2} &0\\0 &0&\begin{array}{cc} 0 & I_{N_3} \\ I_{N_3} &0 \end{array} \end{array}\right).$$
Put $N_1 = \sum_{i=1}^{n_1} \alpha_{s_i}$, $N_2 = \sum_{j=1}^{n_2} \alpha_{t_j}$, $N_3=\sum_{k=1}^{n_3} \alpha_{u_k}$, and $N=N_1+N_2+N_3$, and consider the decomposition of $k^N$ into pairwise orthogonal subspaces $k^{N_1} \oplus k^{N_2} \oplus k^{2N_3}$ given by $\Phi$. Note that representations of $\R$ of dimension $\alpha$ commuting with $\tau$ and $\iota$ adapted to the previous decomposition correspond bijectively to representations of $R$ of the same dimension adapted to this decomposition. This allows us to identify $R(Q,\alpha)$ with the subspace of $R(\Q,\alpha)$ consisting of all representations which preserve the preceding antimorphisms, that is to a subspace of $\X_{\R,N}(k)$, where $\X_{\cdot ,N}$ is the functor introduced in \ref{foncteur} (once the bilinear form has been replaced by $\Phi$).
\subsubsection{} Consider now the algebra $\St_n$ defined as the quotient of $k[e_{s_i},e_{t_j},e_{u_k},e_{u_k^\ast}]$ by the ideal $J$ generated by the relations $e_v^2=e_v$, $e_v e_{v'}=0$ if $v \neq v'$, $\sum_v e_v=1$. This algebra is contained in $\R$. The restriction of $\tau$ to this algebra is exactly the antiinvolution described in \ref{associatequiver}, and we have a fairly nice description of $\X_{\St_n,N}$:
$$\X_{\St_n,N}=\bigcup_{\sigma,\omega} \X_{\sigma,\omega}\, ,$$ 
\noindent where $\sigma$ and $\omega$ range over pairs of admissible vectors in $\mathbb N^n$ such that $\sum \sigma_j=N_1+2N_3$ and $ \sum \omega_j=N_2$, the component $\X_{\sigma,\omega}$ being isomorphic to
$$\bigl(\o_{N_1+2N_3}\times\Sp_{N_2}\bigr) / \Bigl(\prod \bigl(\o_{\sigma_{s_i}} \times \Sp_{\omega_{s_i}}\bigr) \times \prod \bigl(\o_{\sigma_{t_j}} \times \Sp_{\omega_{t_j}}\bigr) \times \prod \bigl(\gl_{\sigma_{u_k}} \times \gl_{\omega_{u_k}}\bigr)\Bigr).$$

It induces a decomposition of $\X_{\R,N}$ as the union $\bigcup_{\sigma,\omega} \varpi^{-1}\X_{\sigma,\omega}$, where we call $\varpi \colon \X_{\R,N} \to \X_{\St_n,N}$ the map induced by the inclusion $\St_n \subset \R$. By applying the argument of \cite[\S3]{LBP} to the component corresponding to the dimensional vectors $\sigma$ and $\omega$ whose coordinate are $\sigma_{s_i}=\alpha_{s_i}$, $\sigma_{t_j}=0$, $\sigma_{u_k}=\alpha_{u_k}$, $\omega_{t_j}=\alpha_{t_j}$ and $\omega_{s_i}=\omega_{u_k}=0$, one gets the expected assertion:

\begin{stheoreme} \label{quiver}
The algebra of polynomials on $R(Q,\alpha)$ invariant under the action of $\prod \o_{\alpha_{s_i}} \times \prod \Sp_{\alpha_{t_j}} \times \prod \gl_{\alpha_{u_k}}$ is generated by the functions $$(f_a)_{{}_a} \mapsto \tr (f_{{\tilde a}_p} \cdots f_{{\tilde a}_1}),$$
\noindent where ${\tilde a}_i$ is an arrow in the associated quiver $\Q$ equals to either $a_i$ or ${a_i}^\ast$, in such a way that $({\tilde a}_1,\ldots,{\tilde a}_p)$ forms an oriented path in that quiver and $f_{{\tilde a}_i}$ means $f_{a_i}$ or its adjoint according to whether ${\tilde a}_i$ is $a_i$ or ${a_i}^\ast$.
\end{stheoreme}

\begin{srem}  It is now easy to deal with the case where we let the whole linear group act (by conjugation) above some of the unpaired vertices: let us call $r_1,\ldots,r_{n_4}$ this vertices, $Q'$ the quiver obtained by adding $n_4$ new vertices ${r}_l^\ast$ and $\alpha ' \in \mathbb N^{n_1+n_2+2(n_3+n_4)}$ the admissible vector naturally deduced from any given admissible dimension vector $\alpha \in \mathbb N^{n_1+n_2+2n_3+n_4}$. The group 
$$\Gamma_\alpha=\prod \o_{\alpha_{s_i}} \times \prod \Sp_{\alpha_{t_j}} \times \prod \gl_{\alpha_{u_k}} \times \prod \gl_{\alpha_{r_l}}$$
\noindent acts on $R(Q,\alpha)$ and $R(Q',\alpha ')$ (the action of $g \in \gl_{\alpha_{r_l}}$ on ${r}_l^\ast$ being $f \mapsto \transp g^{-1} f \transp g$), and the $\Gamma_\alpha$-equivariant projection $k[R(Q',\alpha')] \to k[R(Q,\alpha)]$ allows us to compute the ring of $\Gamma_\alpha$-invariants of $k[R(Q,\alpha)]$.
\end{srem}

\section{Local study of the forgetful map}

In order to simplify the local study of $\mso \to \msl$ it is convenient to investigate separately the injective morphism $\mo \to \mgl$ and the natural map from $\mso$ to the subscheme $\mo^{\O} \subset \mo$ consisting of all orthogonal bundles with trivial determinant. This distinction seems to be quite valuable since the direct differential study of $\mso$ would involve invariant theory for special orthogonal groups, which is far more difficult to deal with (see \ref{gitforso}). We show here that the former is an embedding, while the later is an isomorphism (resp. a $2$-sheeted cover) when $r$ is odd (resp. even).

\subsection{Differential behaviour of $\mo \to \mgl$}

\subsubsection{}\label{luna} Let us now briefly point out the classical way to analyse the local behaviour of $\mo \to \mgl$. Recall first that this application arises as a quotient (by a general linear group $\Gamma=\gl_M$) of an equivariant map between two well-known parameters spaces $R_{\o_r} \to R_{\glr}$ (cf. \cite[7.3]{BLS}). Luna's \'etale slice theorem and deformation theory then allow us to grasp the local structure of these good quotients (cf. \cite[2.5]{kls}): at any polystable vector bundle $E$, $\mgl$ is \'etale locally isomorphic to an \'etale neighbourhood of the origin in the good quotient 
$$\mathrm{Ext}^1(E,E) \git \mathrm{Aut}_{\glr}(E),$$
while $\mo$ is, at any polystable orthogonal bundle $P$, \'etale locally isomorphic to an \'etale neighbourhood of the origin in 
$$H^1(C,\mathrm{Ad}(P)) \git \mathrm{Aut}_{\o_r}(P),$$
where $\mathrm{Ad}(P)$ stands for the vector bundle $P \times^{\o_r} \sok_r$ associated to the adjoint representation of $\o_r$, which is nothing else than the vector bundle of germs of endomorphisms $f$ of $E$ such that  $\sigma f + f^\ast \sigma=0$, where $\sigma \colon E \to E^\ast$ is the symmetric isomorphism given by the quadratic structure on $E$; in other words the adjoint vector bundle $\mathrm{Ad}(P)$ is canonically isomorphic to $\mathbf{\Lambda}^2E^\ast$.

Then, if $P \in \mo$ is a polystable orthogonal bundle with associated vector bundle $E \in \mgl$, the application $\mo \to \mgl$ coincides at $P$, through the preceding local isomorphisms (in the \'etale topology), with the natural map $$H^1(C,\mathrm{Ad}(P)) \git \mathrm{Aut}_{\o_r}(P) \to \mathrm{Ext}^1(E,E) \git \mathrm{Aut}_{\glr}(E)$$ at the origin. In particular the corresponding tangent maps are identified.

\subsubsection{}\label{immersionorthogonal} A more explicit description of the vector spaces $H^1(C,\mathrm{Ad}(P))$ and $\mathrm{Ext}^1(E,E)$ is then strongly needed in order to understand their quotients; we show here that $H^1(C,\mathrm{Ad}(P)) \git \mathrm{Aut}_{\o_r}(P)$ is a closed subscheme of $\mathrm{Ext}^1(E,E) \git \mathrm{Aut}_{\glr}(E)$, which implies our main theorem. 

According to \ref{pointsfermes} the orthogonal structure on the polystable vector bundle $E$ associated to any point $q \in \mo$ gives rise to a splitting of $E$ as a direct orthogonal sum of the form
\begin{align*} \label{polystable} \snum E=\bigoplus\limits_{i=1}^{n_1} E^{(1)}_i \oplus \bigoplus\limits_{j=1}^{n_2}E^{(2)}_j \oplus \bigoplus\limits_{k=1}^{n_3}E^{(3)}_k,
\end{align*}
\noindent where each direct summand $E^{(a)}_l$ may be written as
\begin{itemize}
\item $E^{(1)}_i={F^{(1)}_{i}}\otimes V^{(1)}_{i}$, where $(F^{(1)}_{i})^{}_i$ are mutually non isomorphic stable orthogonal bundles and $(V^{(1)}_{i})^{}_i$ some quadratic vector spaces,
\item $E^{(2)}_j={F^{(2)}_{j}}\otimes V^{(2)}_{j}$, where $(F^{(2)}_{j})^{}_j$ are mutually non isomorphic stable symplectic bundles and $(V^{(2)}_{j})^{}_j$ some symplectic vector spaces,
\item $E^{(3)}_k=({F^{(3)}_{k}} \oplus {F^{(3)}_{k}}{}^\ast) \otimes V^{(3)}_{k}$, where $(F^{(3)}_{k})^{}_k$ are mutually non isomorphic stable bundles such that $F^{(3)}_{k} \not\simeq {F^{(3)}_{k'}}{}^\ast$ and $(V^{(3)}_{k})^{}_k$ some vector spaces carrying a non degenerate bilinear form.
\end{itemize}
\noindent Let us denote by $\psi^{(a)}_{l} \colon F^{(a)}_{l} \to {F^{(a)}_l}{}^\ast$ the duality isomorphism (when $a=1,2$), and $\sigma_l^{(a)} \colon E_{l}^{(a)} \to {E_{l}^{(a)}}{}^\ast$ the symmetric isomorphism defined on $E^{(a)}_l$ (note that ${F^{(3)}_{k}} \oplus {F^{(3)}_{k}}{}^\ast$ has been tacitly endowed with the hyperbolic form). 

The two isotropy groups are then easily identified: we find that 
$$\mathrm{Aut}_{\glr}(E) \simeq \prod\limits_{i=1}^{n_1} \gl({V^{(1)}_{i}}) \times \prod\limits_{j=1}^{n_2} \gl({V^{(2)}_{j}}) \times \prod\limits_{k=1}^{n_3}\Bigl(\gl(V^{(3)}_{k}) \times \gl({V^{(3)}_{k}})\Bigr),$$ 
\noindent while $\mathrm{Aut}_{\o_r}(P) \subset \mathrm{Aut}_{\glr}(E)$ is isomorphic to the subgroup 
$$\prod\limits_{i=1}^{n_1} {\o}({V^{(1)}_{i}}) \times \prod\limits_{j=1}^{n_2} \Sp(V^{(2)}_{j}) \times \prod\limits_{k=1}^{n_3} \gl(V^{(3)}_{k}),$$ 
\noindent where $\gl({V^{(3)}_{k}})$ stands for its image in $\gl({V^{(3)}_{k}}) \times \gl({V^{(3)}_{k}})$ by the morphism $g \mapsto (g, \transp g^{-1})$.

The space $\mathrm{Ext}^1(E,E)$ splits into a direct sum of the spaces $\mathrm{Ext}^1(E^{(k)}_i,E^{(l)}_j)$, and each of these summands is isomorphic to $\mathrm{Ext}^1(F^{(k)}_{i},F^{(l)}_{j}) \otimes \mathrm{Hom}({V^{(k)}_{i}},{V^{(l)}_{j}})$ when neither $k$ nor $l$ equals $3$, or to a sum of summands of this form otherwise. The isotropy groups act on each of those spaces via the natural actions of $\gl(V) \times \gl(V')$ on $\mathrm{Hom}(V,V')$.

An element $\omega=\sum \omega^{(k,l)}_{i,j} \in \mathrm{Ext}^1(E,E) \simeq \bigoplus \mathrm{Ext}^1(E^{(k)}_i,E^{(l)}_j)$ belongs to the space $H^1(C,\mathrm{Ad}(P))$ if and only if $\omega^{(i,i)}_{k,k} \in H^1(C,\boldsymbol{\Lambda}^2E^{(i)}_k{}^\ast) \subset \mathrm{Ext}^1(E^{(i)}_k,E^{(i)}_k)$ and, for $(i,k) \neq (j,l)$, $\sigma^{(l)}_j \omega^{(k,l)}_{i,j} + {\omega^{(l,k)}_{j,i}}{}^\ast \sigma^{(k)}_i=0$. So, identifying $\mathrm{Ext}^1(E_i^{(k)},E_j^{(l)})$ with its image in $\mathrm{Ext}^1(E_i^{(k)},E_j^{(l)}) \oplus \mathrm{Ext}^1(E_j^{(l)},E_i^{(k)})$ by the application $\omega_{i,j}^{(k,l)} \mapsto \omega_{i,j}^{(k.l)} - {\sigma_i^{(k)}}{}^{-1}{\omega_{i,j}^{(k,l)}}{}^\ast \sigma_j^{(l)}$, it appears that $H^1(C,\mathrm{Ad}(P))$ is the subspace of $\mathrm{Ext}^1(E,E)$ isomorphic to the direct sum \stepcounter{seq}
\begin{multline*} 
\snum \label{H1} \bigoplus\limits_{k} \left(\bigoplus\limits_i  H^1(C,\boldsymbol{\Lambda}^2E^{(i)}_k{}^\ast) \oplus \bigoplus\limits_{i<j} \mathrm{Ext}^1(E^{(k)}_i,E^{(k)}_j) \right) \oplus 
\\
\bigoplus\limits_{k<l} \bigoplus\limits_{i,j}  \mathrm{Ext}^1(E^{(k)}_i,E^{(l)}_j),
\end{multline*}
\noindent each one of the diagonal summand being more precisely expressed as:
\begin{multline*} \snum   H^1(C,\boldsymbol{\Lambda}^2E^{(1)}_i{}^\ast)= \Bigl( H^1(C,{\mathbf{S}}^2F^{(1)}_i{}^\ast) \otimes \sok({V^{(1)}_{i}})\Bigr)  \oplus\\ \Bigl(H^1(C,\boldsymbol{\Lambda}^2F^{(1)}_i{}^\ast) \otimes \mathbf{S}^2 {V^{(1)}_{i}}{}^\ast\Bigr),
\end{multline*}
\begin{multline*}
\snum  H^1(C,\boldsymbol{\Lambda}^2F^{(2)}_j{}^\ast)=\Bigl( H^1(C,\boldsymbol{\Lambda}^2F^{(2)}_j{}^\ast) \otimes \spk({V^{(2)}_{j}})\Bigr) \oplus \\ 
\Bigl( H^1(C,\mathbf{S}^2F^{(2)}_j{}^\ast) \otimes \boldsymbol{\Lambda}\!^2 {V^{(2)}_{j}}{}^\ast\Bigr),
\end{multline*}
\begin{multline*}
\snum \label{contributionhyperbolique}  H^1(C,\boldsymbol{\Lambda}^2E^{(3)}_k{}^\ast)= \Bigl(\mathrm{Ext}^1(F^{(3)}_{k},F^{(3)}_{k}) \otimes \glk({V^{(3)}_{k}})\Bigr) \oplus \\ \Bigl(H^1(C,\mathbf{S}^2F^{(3)}_k{}^\ast) \otimes \boldsymbol{\Lambda}^2V^{(3)}_{k}{}^\ast\Bigr) \oplus \Bigl(H^1(C,\boldsymbol{\Lambda}^2F^{(3)}_k{}^\ast) \otimes \mathbf{S}^2({V^{(3)}_{k}}{}^\ast)\Bigr) \oplus \\ \Bigl(H^1(C,\mathbf{S}^2F^{(3)}_k) \otimes \boldsymbol{\Lambda}^2V^{(3)}_{k}{}^\ast\Bigr) \oplus \Bigl(H^1(C,\boldsymbol{\Lambda}^2F^{(3)}_k) \otimes \mathbf{S}^2{V^{(3)}_{k}}{}^\ast\Bigr),
\end{multline*}
\noindent where $\mathrm{Ext}^1(F^{(3)}_{k},F^{(3)}_{k})$ has been identified with its image in $\mathrm{Ext}^1(F^{(3)}_{k},F^{(3)}_{k}) \oplus \mathrm{Ext}^1({F^{(3)}_{k}}{}^\ast,{F^{(3)}_{k}}{}^\ast)$ by the map $\omega \mapsto \omega-\omega^\ast$. Note that the dimensions of all the extension spaces under consideration are trivially available. 

\subsubsection{} This pretty intricate situation suitably expresses itself in terms of representations of quivers. Indeed let us consider the quiver $\mathcal Q$ whose set of vertices is
$${{\mathcal Q}_0=\lbrace{s_{1}^{(1)},\ldots,s^{(1)}_{n_1},s^{(2)}_{1},\ldots,s^{(2)}_{n_2},s^{(3)}_{1},s^{(3^\ast)}_{1},\ldots,s^{(3)}_{n_3},s^{(3^\ast)}_{n_3}\rbrace}},$$ 
\noindent these vertices being connected by $\dim \mathrm{Ext}^1(F_{i}^{(k)},F_{j}^{(l)})$ arrows from $s^{(k)}_i$ to $s^{(l)}_j$ (where we have set $F_{i}^{(3^\ast)}={F_{i}^{(3)}}{}^\ast$). Next define $\alpha \in \mathbb{N}^{n_1 + n_2 +2 n_3}$ according to the dimensions of the corresponding vector spaces $V^{(a)}_{l}$. Therefore the $\mathrm{Aut}_{\gl_r}(E)$-module $\mathrm{Ext}^1(E,E)$ is exactly the $\gl(\alpha)$-module $R(\mathcal Q,\alpha)$ composed of all the representations of $\mathcal Q$ of dimension $\alpha$, and the result of \cite{LBP} recalled earlier provides us with a description of the algebra $k[\mathrm{Ext}^1(E,E)]^{\mathrm{Aut}_{\gl_r}(E)}$.

The inclusion $H^1(C,\mathrm{Ad}(P)) \hookrightarrow \mathrm{Ext}^1(E,E)$ is an $\mathrm{Aut}_{\o_r}(P)=\prod\limits_{i=1}^{n_1} {\o}({V^{(1)}_{i}}) \times \prod\limits_{j=1}^{n_2} \Sp(V^{(2)}_{j}) \times \prod\limits_{k=1}^{n_3} \gl(V^{(3)}_{k})$-equivariant application, so that we have an exact sequence $$k[\mathrm{Ext}^1(E,E)]^{\mathrm{Aut}_{\o_r}(P)} \to k[H^1(C,\mathrm{Ad}(P))]^{\mathrm{Aut}_{\o_r}(P)} \to 0.$$
\noindent This sequence and the theorem \ref{quiver} result in a description of a set of generators for the algebra $k[H^1(C,\mathrm{Ad}(P))]^{\mathrm{Aut}_{\o_r}(P)}$, namely the $(f_a)_{{}_a} \mapsto \tr (f_{{\tilde a}_p} \cdots f_{{\tilde a}_1})$, where $f_{{\tilde a}_i}$ stands for either $f_{a_i}$ or its adjoint map. Now, according to (\ref{H1}),  $H^1(C,\mathrm{Ad}(P))$ is a subspace of $R(\mathcal Q,\alpha)$ made up of representations having the following property: if $f_a \colon V_v \to V_{v'}$ denotes the map associated to an arrow $a \colon v \to v'$, then its adjoint map $f_a^\ast \colon V_{v'}^\ast \to V_{v}^\ast$ is, up to the sign, the map associated to one of the arrows from $v'$ to $v$. So the algebra $k[H^1(C,\mathrm{Ad}(P))]^{\mathrm{Aut}_{\o_r}(P)}$ is generated by traces along oriented cycles in the quiver $\mathcal Q$. This exactly means that the application $k[\mathrm{Ext}^1(E,E)]^{\mathrm{Aut}_{\gl_r}(E)} \to k[H^1(C,\mathrm{Ad}(P))]^{\mathrm{Aut}_{\o_r}(P)}$ is onto.

In view of what has been discussed in \ref{luna} this proves the injectivity of the tangent map of $\mo \to \mgl$ at $q=[P]$. But we have shown in \ref{pointsfermes} that the map $\mo(k) \to \mgl(k)$ is injective. This implies the following: 

\begin{stheoreme}
The forgetful map $(E,q) \in \mo \mapsto E \in \mgl$ is an embedding.
\end{stheoreme}
 
 One easily gets in the very same way the corresponding assertion relative to the moduli of symplectic bundles:
\begin{stheoreme}
The forgetful map $\mathcal M_{{\Sp}_{2r}} \to \mathcal M_{\ssl_{2r}}$ is an embedding.
\end{stheoreme}
The point is that any closed point of $\mathcal{M}_{{\Sp}_{2r}}$ represents a polystable vector bundle of the form 
$$E=\bigoplus\limits_i \Bigl(F^{(1)}_{i}\otimes V^{(1)}_{i}\Bigr) \oplus \bigoplus\limits_j \Bigl(F^{(2)}_{j} \otimes V^{(2)}_{j}\Bigr) \oplus \bigoplus\limits_k \Bigl((F^{(3)}_{k} \oplus F_{k}^{(3) \, \ast})\otimes V^{(3)}_{k}\Bigr),$$
\noindent where $(F^{(1)}_{i})_i$ (resp. $(F^{(2)}_{j})_j$, resp. $(F^{(3)}_{k})_k$) is a family of mutually non isomorphic symplectic (resp. orthogonal, resp. not self-dual) bundles (which are stable as vector bundles), and $(V^{(1)}_{i})_i$ (resp. $(V^{(2)}_{j})_j$, resp. $(V^{(3)}_{k})_k$) a family of quadratic (resp. symplectic, resp. endowed with a non-degenerate bilinear pairing) vector spaces ($F^{(3)}_{k} \oplus F^{(3) \, \ast}_{k}$ being now equipped with the standard symplectic form). Let us denote by $\sigma \colon E \to E^\ast$ the resulting symplectic form on $E$. The bundle $\mathrm{Ad}(P)$ is now isomorphic to the bundles of germs of symmetric endomorphisms of $E$  (that is endomorphisms verifying $\sigma f + f^\ast \sigma=0$), and both the space $H^1(C,\mathrm{Ad}(P))$ and the considered isotropy groups can be described in a manner analogous to that of \ref{immersionorthogonal} (one only has to switch the factors $\boldsymbol{\Lambda}^2F^{(a)}_l{}^\ast$ and $\mathbf{S}^2F^{(a)}_l{}^\ast$, and of course to redefine in the obvious way every map of the form $\mathrm{Ext}^1(F,F') \to \mathrm{Ext}^1(F,F') \oplus \mathrm{Ext}^1(F'^\ast,F^\ast)$). The theorem \ref{quiver} then allows us to conclude again.

\subsection{About $\mso \to \mo$} \label{gitforso} 

\subsubsection{} We have recalled in \ref{giraud} how to compute the fibers of the finite morphism from $\mso$ onto $\mo^{{\O}}=\det^{-1}({\O}_C)$. A point $[P] \in \mo$ in its image has two antecedents if and only if $\mathrm{Aut}_\sor(P) \hookrightarrow \mathrm{Aut}_{\o_r}(P)$ is an isomorphism, that is if and only if every orthogonal bundle $F_{i}^{(1)}$ appearing in the splitting (\ref{polystable}) of $E$ has even rank.

Luna's theorem reduces once again the differential study of this application to an invariant calculus: the tangent map of $\mso \to \mo$ at $[P] \in \mso$ is indeed identified with that of  $H^1(C,\mathrm{Ad}(P)) \git \mathrm{Aut}_\sor(P) \to H^1(C,\mathrm{Ad}(P)) \git \mathrm{Aut}_{{\o}_r}(P)$ (at the origin).

Therefore, if $r$ is odd, $\mso \to \mo^{\O}$ is an isomorphism.

\subsubsection{} Let us consider now the even case. The morphism $\mso \to \mo$ is then a $2$-sheeted cover, which is \'etale above the locus of points having two antecedents. A branched point corresponds to an orthogonal polystable bundle $E$ containing at least one subbundle isomorphic to $F^{(1)}_{i} \otimes V^{(1)}_{i}$ where $F^{(1)}_{i}$ is an orthogonal bundle of odd rank: we then have to understand the inclusion 
$$k[H^1(C,\mathrm{Ad}(P))]^{\mathrm{Aut}_{\o_r}(P)} \hookrightarrow k[H^1(C,\mathrm{Ad}(P))]^{\mathrm{Aut}_{\sor}(P)}.$$ 
\noindent It is easy to produce a primitive element for the generic extension (which is of degree $2$). First note that the vector space $W$ obtained as the direct sum of the $V^{(1)}_i$ corresponding to the orthogonal bundles $F^{(1)}_i$ of odd rank has even dimension, and has an orthogonal structure inherited from the ones of the $V^{(1)}_i$. The space composed of all the antisymmetric endomorphisms of $W$ may be identified with a direct summand of $H^1(C,\mathrm{Ad}(P))$, and mapping any element $\omega \in H^1(C,\mathrm{Ad}(P))$ to the pfaffian of the endomorphism of $W$ induced by $\omega$ then defines a function belonging to $k[H^1(C,\mathrm{Ad}(P))]^{\mathrm{Aut}_\sor(P)}$ which is not $\mathrm{Aut}_{\o_r}(P)$-invariant; this function certainly generates the generic extension. 

It is more difficult to give a convenient description of this algebra: in the (simplest) case where $P$ is isomorphic to $\O \otimes V$ with $V$ a quadratic vector space of even dimension, we have to understand the action of $\mathrm{Aut}_{\sor}(P) \simeq \sor$ on $H^1(C,\mathrm{Ad}(P)) \simeq H^1(C,\O) \otimes \sok(V)$. This can be solved again thanks to Procesi's trick (cf. \ref{procesitrick}): the computation has been carried out in \cite{aslaksen}, and provides a set of generators for the $k[H^1(C,\mathrm{Ad}(P))]^{\mathrm{Aut}_{\o_r}(P)}$-algebra $k[H^1(C,\mathrm{Ad}(P))]^{\mathrm{Aut}_{\sor}(P)}$ in terms of \textit{polarized pfaffians}. 

Let us finally mention that in the general case we can easily infer from the main result of \cite{lopatin} a family of generators for $k[H^1(C,\mathrm{Ad}(P))]^{\mathrm{Aut}_\sor(P)}$ which are also obtained as polarized pfaffians.

\subsection{About the multiplicity at stable points}

\subsubsection{} The discussion held in \ref{luna} contains in fact a more precise statement, related to the completed local rings of $\mo$ and $\mgl$. Indeed, if $q$ is a point of $\mo$ representing a polystable bundle $P$ whose image in $\mgl$ is a point $s=[E]$, we have the following commutative diagram, where the rings of the second row are the completions of the local rings (of the involved algebras of invariants) at the origin, 
\begin{align} \snum  \label{lunadiag} 
\begin{array}{c}
\xymatrix@C=35pt{
\widehat{\O}_{\mgl ,\, s} \ar[d]^-{\wr} \ar[r]  &\widehat{\O}_{\mo ,\, q} \ar[d]^-\wr\\
 \Bigl(k[\mathrm{Ext}^1(E,E)]^{\mathrm{Aut}_{\gl_r}(E)}\Bigr)^{\widehat{\vphantom{}}}\ar[r] & \Bigl(k[H^1(C,\mathrm{Ad}(P))]^{\mathrm{Aut}_{\o_r}(P)}\Bigr)^{\widehat{\vphantom{}}}.}
\end{array}
\end{align}
\noindent This description of the completed local rings of $\mo$ provides us with additional informations about the local structure of $\mo$, at least at the points where the situation is not too bad (see \cite{localstr} for the case of $\mgl$): the more we know about the second main theorem of invariant theory for the isotropy group of $P$, the easier our calculations will be.

\subsubsection{} Let $P$ be an orthogonal bundle whose underlying vector bundle is of the form $E=E_1 \oplus E_2$, with $E_1$ and $E_2$ two non-isomorphic $\gl$-stable orthogonal bundles. The description of the inclusion $H^1(C,\mathrm{Ad}(P)) \hookrightarrow \mathrm{Ext}^1(E,E)$ given in (\ref{H1}) here reduces to
$$\begin{array}{ccccccc}
H^1(C,\mathrm{Ad}(P)) &=& H^1(C,\boldsymbol{\Lambda}^2E_1{}^\ast) & \oplus & \mathrm{Ext}^1(E_1,E_2) & \oplus & H^1(C,\boldsymbol{\Lambda}^2E_2{}^\ast) \\
\cap                  & &    \cap                 &        &          \cap           &        &   \cap     \\
\mathrm{Ext}^1(E,E)   &=& \mathrm{Ext}^1(E_1,E_1) & \oplus & \mathrm{Ext}^1(E_1,E_2)\  \oplus \ \mathrm{Ext}^1(E_2,E_1) & \oplus & \mathrm{Ext}^1(E_2,E_2); \end{array}$$
\noindent the isotropy subgroup $\mathrm{Aut}_{\o_r}(P)$, isomorphic to $\mathbb Z / 2 \mathbb Z \times \mathbb Z/2\mathbb Z$, acts trivially on $H^1(C,\boldsymbol{\Lambda}^2E_i{}^\ast)$ and by multiplication by $\pm 1$ on $\mathrm{Ext}^1(E_1,E_2)$ (while $(\alpha_1,\alpha_2) \in \mathrm{Aut}_\glr(E) \simeq \mathbb G_m \times \mathbb G_m$ acts on $\mathrm{Ext}^1(E_i,E_j)$ by $\alpha_j \alpha_i^{-1}$).

Let $(X_k^{(i)})^{}_k$ (resp. $(Y_l)_l$) be a basis of $ H^1(C,\boldsymbol{\Lambda}^2E_i{}^\ast)^\ast$ (resp. $\mathrm{Ext}^1(E_1,E_2)^\ast \subset H^1(C,\mathrm{Ad}(P))^\ast$). Then $k[H^1(C,\mathrm{Ad}(P))]^{\mathrm{Aut}_{\o_r}(P)}$ is the subring of $k[X_k^{(i)},Y_l]$ generated by all the $X_k^{(i)}$ and the products $Y_l Y_{l'}$; if $\mathcal V$ denotes the affine cone over the Veronese variety $\mathbb P(\mathrm{Ext}^1(E_1,E_2)) \subset \mathbb P\left(\mathbf{S}^2\mathrm{Ext}^1(E_1,E_2)\right)$ we get the following isomorphism:
$$\mathrm{Spec}(k[H^1(C,\mathrm{Ad}(P))]^{\mathrm{Aut}_{\o_r}(P)}) \buildrel\sim\over\lra \Bigl(H^1(C,\mathrm{Ad}(P_1))\oplus H^1(C,\mathrm{Ad}(P_2))\Bigr) \times \mathcal V.$$

Using the identification $\widehat{\O}_{\mo,\, q} \simeq \Bigl(k[H^1(C,\mathrm{Ad}(P))]^{\mathrm{Aut}_{\o_r}(P)}\Bigr)^{\widehat{\vphantom{}}}$ we have the following result:

\begin{scor}
The tangent space to $\mo$ at a point $[P]$ given as the direct sum $E_1 \oplus E_2$ of two non-isomorphic stable orthogonal bundles is isomorphic to
$$ H^1(C,\boldsymbol{\Lambda}^2E_1{}^\ast)\oplus  H^1(C,\boldsymbol{\Lambda}^2E_2{}^\ast) \oplus (\mathbf{S}^2\mathrm{Ext}^1(E_1,E_2)),$$
and the multiplicity of $\mo$ at this point is equal to $\displaystyle{2 ^{r_{1} r_{2}(g-1)-1}}$, where $r_i$ is the rank of $E_i$.
\end{scor}

\begin{srem} (i) The general case of a stable point $q \in \mo$ is more difficult: such a bundle corresponds to a vector bundle which splits as a direct sum of $n$ mutually non-isomorphic $\gl_{r_i}$-stable orthogonal vector bundles. We can use \ref{quiver} to try to get some additional informations about the local structure at $q$, for instance by computing the multiplicity at this point. One can easily check that, if $n=3$ (resp. $4$) this multiplicity is equal to $\displaystyle{2\prod_{i<j}2^{r_ir_j(g-1)-1}}$ (resp. $\displaystyle{8\prod_{i<j}2^{r_ir_j(g-1)-1}}$).

\noindent(ii) It is not hard to deal with the case of an orthogonal (non stable) bundle of the form $F \oplus F^\ast$ with $F \not\simeq F^\ast$: we see that $\mathcal M_{\o_r}$ is, at such a point, \'etale locally isomorphic to $\mathrm{Ext}^1(F,F) \oplus \mathcal S$, where $\mathcal S$ is the affine cone over the Segre variety 
$$\mathbb P(H^1(C, \boldsymbol{\Lambda}^2F^\ast))\times \mathbb P(H^1(C,\boldsymbol{\Lambda}^2F))\subset \mathbb P\left(H^1(C,\boldsymbol{\Lambda}^2F^\ast)\otimes H^1(C,\boldsymbol{\Lambda}^2F)\right).$$
\end{srem}

\nocite{*}
\bibliographystyle{amsplain}
\bibliography{immersion}

\end{document}